\documentclass[letterpaper,11pt]{amsart}
\usepackage{amsmath,amssymb,amsthm,pinlabel,tikz,hyperref,mathrsfs,thmtools}
\usepackage{verbatim}
\usepackage{float}

\usetikzlibrary{patterns}

\newcommand{\nc}{\newcommand}
\nc{\dmo}{\DeclareMathOperator}
\dmo{\ra}{\rightarrow}
\dmo{\N}{\mathbb{N}}
\dmo{\Z}{\mathbb{Z}}
\dmo{\R}{\mathbb{R}}
\dmo{\T}{\mathcal{T}}
\dmo{\C}{\mathcal{C}}
\dmo{\AC}{\mathcal{AC}}
\dmo{\Mod}{Mod}
\dmo{\PMod}{PMod}
\dmo{\B}{B}
\dmo{\PB}{PB}
\dmo{\Sp}{Sp}
\dmo{\I}{\mathcal{I}}
\dmo{\el}{\ell_{\C}}
\dmo{\NN}{\mathcal{N}}
\dmo{\Tr}{Tr}
\dmo{\llangle}{\langle\langle}
\dmo{\rrangle}{\rangle\rangle}

\usepackage{tikz-cd}

\tikzcdset{row sep/normal=1.5cm}
\tikzcdset{column sep/normal=1.5cm}

\usepackage{caption}
\usepackage{subcaption}

\usetikzlibrary{decorations.markings}
\tikzset{->-/.style={decoration={
			markings,
			mark=at position #1 with {\arrow{>}}},postaction={decorate}}}

\nc{\nt}{\newtheorem}

\nt{theorem}{Theorem}

\newtheorem{thm}{{\bf Theorem}}[section]
\newtheorem{lem}[thm]{{\bf Lemma}}
\newtheorem{cor}[thm]{{\bf Corollary}}
\newtheorem{prop}[thm]{{\bf Proposition}}

\newtheorem{remark}[thm]{Remark}

\newtheorem{question}[thm]{Question}
\newtheorem{conj}[thm]{Conjecture}

\numberwithin{equation}{section}

\title[Minimal asymptotic
  translation lengths on curve complexes]{Minimal asymptotic
  translation lengths\\on curve complexes\\and homology of mapping~tori}

\date{\today}
\author{Hyungryul Baik}
\address{
	Department of Mathematical Sciences, KAIST\\
	291 Daehak-ro Yuseong-gu, Daejeon, 34141, South Korea 
}
\email{%
	hrbaik@kaist.ac.kr
}

\author{Dongryul M. Kim}

\address{%
	Department of Mathematics, Yale University\\
	219 Prospect Street, New Haven, CT 06511, USA
}
\email{%
	dongryul.kim@yale.edu
}

\author{Chenxi Wu}

\address{%
	Department of Mathematics, University of Wisconsin--Madison\\
	480 Lincoln Drive, Madison, WI 53706, USA
}
\email{%
	cwu367@math.wisc.edu
}

\begin{document}
\begin{abstract}
	Let $S_g$ be a closed orientable surface of genus $g > 1$. Consider the minimal asymptotic translation length $L_{\mathcal{T}}(k, g)$ on the Teichm\"uller space of $S_g$, among pseudo-Anosov mapping classes of $S_g$ acting trivially on a $k$-dimensional subspace of $H_1(S_g)$, $0 \le k \le 2g$. The asymptotics of $L_{\mathcal{T}}(k, g)$ for extreme cases $k = 0, 2g$ have been shown by several authors. Jordan Ellenberg asked whether there is a lower bound for $L_{\mathcal{T}}(k, g)$ interpolating the known results on $L_{\mathcal{T}}(0, g)$ and $L_{\mathcal{T}}(2g, g)$, which was affirmatively answered by Agol, Leininger, and Margalit.

	In this paper, we study an analogue of Ellenberg's question, replacing Teichm\"uller spaces with curve complexes. We provide lower and upper bound on the minimal asymptotic translation length $L_{\mathcal{C}}(k, g)$ on the curve complex, whose lower bound interpolates the known results on $L_{\mathcal{C}}(0, g)$ and $L_{\mathcal{C}}(2g, g)$.
	
	Finally, for each $g$, we construct a non-Torelli pseudo-Anosov $f_g \in \Mod(S_g)$ which does not normally generates $\Mod(S_g)$ and so that the asymptotic translation length of $f_g$ on curve complexes decays more quickly than a constant multiple of $1/g$ as $g \to \infty$. From this, we provide a restriction on how small the asymptotic translation lengths on curve complexes should be if the similar phenomenon as in the work of Lanier and Margalit on Teichm\"uller spaces holds for curve complexes.

	\medskip
	{\noindent \bf MSC 2020 classes: 37E30, 57K20}
\end{abstract}

\maketitle

%
%

\section{Introduction}	\label{section:introduction}

Let $S_g$ be a closed connected orientable surface of genus $g > 1$, $\Mod(S_g)$ be its mapping class group, and $\C(S_g)$ be its curve complex. Then $\Mod(S_g)$ isometrically acts on $\C(S_g)$, hence the \emph{asymptotic translation length} $\ell_{\C}(f)$ of $f \in \Mod(S_g)$ on $\C(S_g)$ is defined: $$\ell_{\C}(f) := \liminf_{n\to\infty} {d_{\C}(x, f^n(x)) \over n}$$ for any $x \in \C(S_g)$ where $d_{\C}$ is the standard metric on $\C(S_g)$. The asymptotic translation length is also called \emph{stable translation length}.

Note that $\Mod(S_g)$ also acts on $H_1(S_g)$, the first homology group of $S_g$ with real coefficients. For $f \in \Mod(S_g)$, we denote the dimension of a maximal subspace of $H_1(S_g)$ on which $f$ is trivial by $m(f)$. In particular, $m(f) = 2g$ if and only if $f$ is in the \emph{Torelli group} $\I_g \le \Mod(S_g)$, the subgroup consisting of elements which act trivially on $H_1(S_g)$. As an application of Mayer-Vietoris sequence, one can observe that $m(f) + 1$ is same as the first betti number of the mapping torus of $f$, which is hyperbolic if and only if $f$ is pseudo-Anosov by Thurston \cite{thurston1998hyperbolic}.

In this paper, we mainly study the \emph{minimal asymptotic translation lengths} among pseudo-Anosov mapping classes acting trivially on some subspaces of homology groups. Namely, for $0 \le k \le 2g$, we define $$L_{\C}(k, g) := \inf \left\lbrace \ell_{\C}(f) : f \in \Mod(S_g), f \mbox{ is pseudo-Anosov}, m(f) \ge k \right\rbrace.$$ Then we investigate asymptotics of $L_{\C}(k, g)$ with varying $k$ and $g$.

By replacing the curve complex $\C(S_g)$ with Teichm\"uller space $\T(S_g)$, one can also define $L_{\T}(k, g)$ analogously. Note that $\ell_{\T}(f)$ for a pseudo-Anosov element $f$ is same as the logarithm of the stretch factor \cite{bers1978extremal}, hence coincides with the topological entropy of $f$ \cite[Expos\'e Ten]{fathi2012thurston}.

 In each setting, there are two extreme cases: the first extreme is the case $k = 0$ that the minimal asymptotic translation length is considered in the \emph{entire mapping class group $\Mod(S_g)$.} The other extreme is $k = 2g$, which means that the minimal asymptotic translation length is considered in the \emph{Torelli subgroup $\I_g < \Mod(S_g)$.} These four cases have been resolved by various authors as below:\footnote{Throughout the paper, we write $A(x) \gtrsim B(x)$ if there exists a uniform constant $C > 0$ such that $ A(x) \le CB(x)$ for all $x$ in the domain. We also write $A(x) \asymp B(x)$ if $A(x) \gtrsim B(x)$ and $B(x) \gtrsim A(x)$.}

 \begin{table}[H]
 \begin{center}
 \begin{tabular}{ |c|c|c| } 
	\hline
	& Teichm\"uller spaces & Curve complexes\\ 
	\hline
	$\Mod(S_g)$ &  $\begin{matrix} \\ \mbox{(Penner \cite{penner1991bounds})} \\ L_{\T}(0, g) \asymp 1/g  \\ \ \end{matrix}$ & $\begin{matrix} \mbox{(Gadre-Tsai \cite{GadreTsai11})} \\ L_{\C}(0, g) \asymp 1/g^2 \end{matrix}$ \\ 
	\hline
	$\I_g$  & $\begin{matrix} \\ \mbox{(Farb-Leininger-Margalit \cite{farb2008lower})} \\ L_{\T}(2g, g) \asymp 1 \\ \ \end{matrix}$ & $\begin{matrix} \mbox{(Baik-Shin \cite{baik2020minimal})} \\ L_{\C}(2g, g) \asymp 1/g \end{matrix}$ \\ 
	\hline
   \end{tabular}
   \caption{Four extreme cases of $L_{\T}(k, g)$ and $L_{\C}(k, g)$.} \label{table.torelliandmcg}
\end{center}
\end{table}

Ellenberg \cite{ellenbergpuzzle} asked if $L_{\T}(k, g)$ interpolates $L_{\T}(0, g)$ and $L_{\T}(2g, g)$ in the sense that there exists $C > 0$ such that \begin{equation} \label{eqn:ell}
	L_{\T}(k, g) \geq C(k+1)/g
\end{equation} for all $g > 1$ and $0 \le k \le 2g$. This was answered affirmatively by Agol-Leininger-Margalit in \cite{agol2016pseudo}. Indeed, they actually showed $L_{\T}(k, g) \asymp (k+1)/g$.

We ask an analogous question whether $L_{\C}(k, g)$ interpolates $L_{\C}(0, g)$ and $L_{\C}(2g, g)$ in a similar sense as Ellenberg's question \eqref{eqn:ell}. We show that this is indeed the case and more concretely we obtain the following: 
\begin{thm} \label{thm:main}
	There exists $C, C'>0$ such that $${C\over g(2g-k+1)} \le L_{\C}(k, g) \le C' {k + 1 \over g \log g}$$ for all $g > 1$ and $0 \le k \le 2g$.
\end{thm}

From the statement, if $k$ grows at least $2g - C'$ for some constant $C' > 0$, then $L_{\C}(k, g) \gtrsim 1/g$ while $L_{\C}(0, g) \asymp 1/g^2$. Observing this, we ask about minimal $k$ with $L_{\C}(k, g) \asymp 1/g$. For this discussion, see Section \ref{sec:example}.

While the lower bound in Theorem \ref{thm:main} interpolates $L_{\C}(0, g) \asymp 1/g^2$ and $L_{\C}(2g, g) \asymp 1/g$, the upper bound in Theorem \ref{thm:main} does not interpolate these two values well. Indeed, we construct some values of $k$ and $g$ showing that ${k + 1 \over g \log g}$ is larger than the actual asymptotic. We also show that $k/g^2$ works as an upper bound for some choices of $(k, g)$, which interpolates $L_{\C}(0, g) \asymp 1/g^2$ and $L_{\C}(2g, g) \asymp 1/g$.

\begin{thm} \label{thm:partialupper}
	There is a uniform constant $C>0$ satisfying the following:
	
	For any integers $g, k \ge 0$, there exists a pseudo-Anosov $f : S_{g'} \to S_{g'}$ such that $g' > g$, $m(f) = k' > k$, and $$\ell_{\C}(f) \le C{ k' \over g'^2}.$$
\end{thm}

Applying Theorem \ref{thm:partialupper} inductively, it follows that there is a diverging sequence $(k_j, g_j) \to \infty$ so that $L_{\C}(k_j, g_j) \lesssim k_j/g_j^2$. See Corollary \ref{cor:uppersubseq}. Based on Table \ref{table.torelliandmcg}, we conjecture that the upper bound in Theorem \ref{thm:partialupper} is actually the asymptotic for $L_{\C}(k, g)$.

\begin{conj} \label{conjecture} We have
	$$L_{\C}(k, g) \asymp {k \over g^2}$$ for $g > 1$ and $0 \le k \le 2g$.
\end{conj}

We focus on specific dimensions of maximal invariant subspaces. In \cite{baik2020minimal}, Torelli pseudo-Anosovs are constructed in a concrete way based on Penner's or Thurston's construction. In a similar line of thought, we utilize finite cyclic covers of $S_2$ so that we get pseudo-Anosovs $f \in \Mod(S_g)$ with $m(f) = 2g-1$ and satisfying the upper bound in Theorem \ref{thm:partialupper}. As a consequence, this yields the asymptotic of $L_{\C}(2g-1, g)$; only two extreme cases $\Mod(S_g)$ and $\I_g$ were previously known. It is also interesting to figure out the asymptote $L_{\C}(k, g)$ for other values $(k, g)$:

\begin{question}
	Can we give a sequence $(k_j, g_j)$, other than $(0, g)$ and $(2g, g)$, with explicit asymptote for $L_{\C}(k_j, g_j)$ as $j \to \infty$?
\end{question}

We give one such example in the following:

\begin{thm} \label{thm:2g+1}
	There exists a uniform constant $C > 0$ and pseudo-Anosovs $f_g \in \Mod(S_g)$ such that $$m(f_g) = 2g - 1 \quad \mbox{and} \quad \ell_{\C}(f_g) \le {C \over g}$$ for all $g > 1$. Moreover, the following asymptotic holds: $$L_{\C}(2g-1, g) \asymp {1 \over g}.$$
\end{thm}

The construction involved in Theorem \ref{thm:2g+1} can be modified to deal with the Torelli case. Such a modification gives an asymptotic for $L_{\C}(2g, g)$ which was already shown by \cite{baik2020minimal} in a different way. See Remark \ref{rmk:torelli}. Further, only the last assertion can also be deduced from Theorem \ref{thm:main} and \cite{baik2020minimal}. See Section \ref{sec:example} for details.

In \cite{lanier2018normal}, Lanier and Margalit showed that a pseudo-Anosov with small asymptotic translation length on the Teichm\"uller space has an entire mapping class group as its normal closure. The first and the third authors, Kin, and Shin, made an analogous question for asymptotic translation lengths on curve complexes in \cite{baik2019asymptotic} (see \cite[Question 1.2]{baik2019asymptotic}). We later show that pseudo-Anosovs $f_g$ constructed in Theorem \ref{thm:2g+1} never normally generate the mapping class groups. As $\ell_{\C}(f_g)$ is concretely estimated in Section \ref{sec:example}, it provides how small the asymptotic translation length should be in order to observe the similar phenomenon as in \cite{lanier2018normal}. In other words, we prove the following.

\begin{thm} \label{prop:res}
	Suppose that there exists a universal constant $C$ so that if a non-Torelli pseudo-Anosov $f \in \Mod(S_g)$ has $\ell_{\C}(f) < C/g$ then $f$ normally generates $\Mod(S_g)$ for large $g$. Then $$C \le 1152.$$
\end{thm}

\subsection*{Organization}In Section \ref{sec:main}, we prove Theorem \ref{thm:main}. Theorem \ref{thm:partialupper} is proved in Section \ref{subsec:partial}. In Section \ref{sec:example}, explicit construction of pseudo-Anosovs realizing the asymptotic for $L_{\C}(2g-1, g)$ is provided, implying Theorem \ref{thm:2g+1}. The discussion on small asymptotic translation lengths on curve complexes and normal generation of mapping class groups is provided in Section \ref{sec:normal}.

\subsection*{Acknowledgements}

The authors greatly appreciate Changsub Kim and Yair N. Minsky for helpful discussions.
We also thank the anonymous referee for helpful comments.
The first author was partially supported by the National Research Foundation of Korea(NRF) grant funded by the Korea government(MSIT) (No. 2020R1C1C1A01006912).

%
%

\section{Proof of Theorem \ref{thm:main}} \label{sec:main}

In this section, we prove Theorem \ref{thm:main}.

\subsection*{Lower bound} \label{sec:lower}

The main idea of showing the lower bound is similar to the one used in the proof in \cite{baik2020minimal} of that $L_{\C}(2g, g) \geq C/g$ for some constant. First note that for any $f: S_g \to S_g$, the Lefschetz number $L(f)$ is $2-\Tr(f_\ast)$ where $\Tr(f_\ast)$ is the trace of the induced map $f_*: H_1(S_g)\rightarrow H_1(S_g)$. 

Let us fix a pseudo-Anosov $f: S_g \to S_g$ whose restriction onto a $k$-dimensional subspace of $H_1(S_g)$ is the identity.

Pick a suitable basis for $H_1(S_g)$, the matrix for $f_*$ can be written as
  \[\left(\begin{array}{cc} I_k & *\\ 0 & M\end{array}\right).\]
 Suppose first that $k > 0$. When $k$ is odd, let $m = 2g-k$ and when $k$ is even, let $m= 2g-(k-1)$. By taking the upper left block to be $I_{k-1}$ in case $k$ is even, one may assume $M$ is a $m \times m$ square matrix with determinant $1$ and $m$ is odd (determinant 1 comes from the fact that $f_*$ is actually a symplectic matrix). 

  Recall that there is a relation between trace and determinant as follows: 
  \begin{lem}[{\cite[Appendix B]{kondratyuk1992superconducting}}] \label{lem:dettrace}
  	 For any $m\times m$ matrix $A$,
  \[(-1)^m \det A=\sum\limits_{\substack{c_1, \dots, c_m \geq 0, \\ c_1 + 2c_2 + \cdots mc_m = m}} 
      \prod_{i=1}^m{1\over c_i !}\left(-{\Tr(A^i)\over i}\right)^{c_i}\]
  \end{lem}

Observe that at least one of the matrices $M$, $M^2$, $\dots$, $M^m$ must have positive trace. Otherwise the right-hand side of the equality in Lemma \ref{lem:dettrace} is always non-negative when we plug in $M$ in the place of $A$ in the lemma. On the other hand, since $\det(M)=1$ and $m$ is always odd by our choice, the left-hand side is $-1$, a contradiction. 

This implies that for some $j$ satisfying $1 \leq j \leq m \leq 2g-k+1$, $\Tr(M^j)$ is positive, i.e., at least 1 since it is an integral matrix. $\Tr(f_*^j)$ is the sum of $\Tr(M^j)$ and the trace of the upper left block which is $2g-m \geq 1$. Therefore, $\Tr(f_*^j)$ is at least $2$ in general. But in fact $2g-m \geq 3$ as long as $k \geq 3$. 

Assume $k \geq 3$. Now we have that $L(f^j) = 2 - \Tr(f_*^j) < 0$ and we can apply a result of Tsai \cite{Tsai09}. 
Then $\ell_{\C} (f^j) \geq C/g$ for some constant $C>0$ and consequently, 
$$\ell_{\C} \geq \dfrac{C}{gj} \geq \dfrac{C}{g(2g-k+1)}.$$

For any fixed $k$, $\ell_{\C} \geq L_{\C}(0, g) \asymp 1/g^2 $. Hence $\ell_{\C} \geq \dfrac{C_k}{g(2g-k+1)}$ for some constant $C_k$. Since the above argument works for any $k \geq 3$, replacing $C$ by $\min\{C, C_0, C_1, C_2\}$, we obtain the lower bound in Theorem \ref{thm:main}.

%
%

\subsection*{Upper bound} \label{sec:upper} \label{subsec:upperboundteich}

We now prove the upper bound provided in Theorem~\ref{thm:main}.

Recall that the Teichm\"uller space $\T(S_g)$ is a space of marked hyperbolic structures on $S_g$ and vertices of the curve complex $\C(S_g)$ are essential simple closed curves on $S_g$. Hence, we can associate each point $x \in \T(S_g)$ with systoles on $S_g$, the shortest simple closed geodesics, in the hyperbolic structure $x$. As systoles are within a uniformly bounded distance in the curve complex, it gives a coarsely well-defined map $\pi_g : \T(S_g) \to \C(S_g)$.

Masur-Misnky studied $\pi_g : \T(S_g) \to \C(S_g)$ in \cite{MasurMinsky99} and showed that $\pi_g$ is coarsely Lipschitz. 

\begin{prop}[$(K_g, D_g)$-coarsely Lipschitz, {\cite{MasurMinsky99}}]
	There exist constants $K_g, D_g > 0$ such that for any $x, y \in \T(S_g)$ we have $$d_{\C}(\pi_g(x), \pi_g(y)) \le K_g d_{\T}(x, y) + D_g$$ where $d_{\T}$ is the Teichm\"uller metric.
\end{prop}

Furthermore, $\pi_g$ is coarsely $\Mod(S_g)$-equivariant in the sense that there exists a constant $A_g$ such that $d_{\C}((\pi_g \circ f)(x), (f \circ \pi_g)(x)) \le A_g$ for any $x \in \T(S_g)$ and $f \in \Mod(S_g)$. Then for $f \in \Mod(S_g)$, $n > 0$, and $x \in \T(S_g)$, we have $$
\begin{aligned}
d_{\C}(\pi_g(x), f^n(\pi_g(x))) & \le d_{\C}(\pi_g(x), \pi_g(f^n(x))) + A_g \\
& \le K_g d_{\T}(x, f^n(x)) + D_g + A_g
\end{aligned}$$ Hence, we now have the comparison between asymptotic translation lengths of $f \in \Mod(S_g)$ measured on $\C(S_g)$ and $\T(S_g)$: $$\ell_{\C}(f) \le K_g \ell_{\T}(f)$$ In particular, we have \begin{equation} \label{eqn:estimatefromLip}
L_{\C}(k, g) \le K_g L_{\T}(k, g).
\end{equation}

Due to the work \cite{agol2016pseudo} of Agol-Leininger-Margalit, we already know the asymptotic of $L_{\T}(k, g)$. Hence, it remains to figure out the asymptotic of $K_g$.

In \cite{GadreHironakaKentLeininger13}, Gadre-Hironaka-Kent-Leininger considered the minimal possible Lipschitz constant $K_g$ which is defined as $$\kappa_g := \inf \left\lbrace K_g \ge 0 : \pi_g \mbox{ is } (K_g, D_g)\mbox{-coarsely Lipschitz for some } D_g > 0 \right\rbrace.$$ Then they showed that $$\kappa_g \asymp {1 \over \log g}.$$ Combining this with \cite{agol2016pseudo} and Inequality \eqref{eqn:estimatefromLip}, we deduce the upper bound in Theorem \ref{thm:main}.

\section{Upper bound interpolates $L_{\C}(0, g)$ and $L_{\C}(2g, g)$} \label{subsec:partial}

The upper bound provided in Theorem \ref{thm:main} does not interpolate $L_{\C}(0, g)$ and $L_{\C}(2g, g)$, and it is not sharp enough as one can see in Section \ref{sec:example}. As stated in Theorem \ref{thm:partialupper}, the upper bound conjectured in Conjecture \ref{conjecture} can be observed along certain sequence $(k_j, g_j) \to \infty$. This section is devoted to proving Theorem \ref{thm:partialupper}.

\subsection*{Proof of Theorem \ref{thm:partialupper}}

Let $f_0$ be a pseudo-Anosov map in the Torelli group of genus $g_0 >~1$. Let $M$ be its mapping torus, $\alpha \in H^1(M)$ be the first cohomology class of $M$ corresponding to $f_0$, $\beta$ be an element in $H^1(M)$ which is restricted to a cohomology class dual to a simple closed curve $\gamma$ on $S_{g_0}$. For large enough $n > g + k$, let $f_{n}$ be the pseudo-Anosov monodromy corresponding to $2^n\alpha+\beta$. Then $f_n$ has the fiber of genus $O(2^n)$, and $\ell_{\C}(f_n)$ is $O(2^{-2n})$. (cf. \cite{baik2018upper})

A way to construct the surface $S_n$ and map $f_n$ corresponding to $2^n\alpha+\beta$ is as follows: Let $\widehat{S}$ be the $\mathbb{Z}$-fold cover corresponding to $\beta$ restricted to $S_{g_0}$, $\widehat{f}$ a lift of $f_0$, $h$ the deck transformation, then with a suitable choice of $\widehat{f}$ we have $S_n=\widehat{S}/(h^{2^n}\widehat{f})$ and $f_n$ is lifted to $h$. Now consider a simple closed curve on a fundamental domain of $\widehat{S}$ which is not homologous to the boundary, such that the homology class $c$ represented by this curve $\gamma$ is preserved by $\widehat{f}$. The existence of such a homology class is due to the construction in Baik-Shin \cite{baik2020minimal}. Then $\sum_{i=0}^{2^n-1}f_n^ic$ is invariant under $f_n$, and for $k<n$, let $c_k=\sum_{i=0}^{2^{n-k}-1}f_n^{i2^k}c$. Now $\mathrm{Span}\{c_k, f_nc_k,\dots f_n^{2^k-1}c_k\}$ is a $2^k$ dimensional invariant subspace of $f_n^{2^k}$. This proves Theorem \ref{thm:partialupper}.

\medskip

Since the constant $C$ in Theorem \ref{thm:partialupper} does not depend on the choice of given $g$ and $k$, we can apply the theorem inductively: At each $j$-th step with $g_j$ and $k_j$, Theorem \ref{thm:partialupper} applied to $g_j$ and $k_j$ gives $g' > g_j$, $k' > k_j$, and a pseudo-Anosov $f_{j+1} : S_{g_j'} \to S_{g_j'}$ with $\ell_{\C}(f_{j+1}) \le C {k_{j}'/ g_{j}'^2}$. Then we set $g_{j+1} := g_j'$ and $k_{j+1} := k_j'$. As a consequence, we obtain the following corollary that interpolates $L_{\C}(0, g)$ \cite{GadreTsai11} and $L_{\C}(2g, g)$ \cite{baik2020minimal} in a partial way.

\begin{cor} \label{cor:uppersubseq}
	There is a constant $C$ and a diverging sequence $(k_j, g_j) \to~\infty$ as $j \to \infty$ such that $$L_{\C}(k_j, g_j) \le C {k_j \over g_j^2}.$$
\end{cor}

Corollary \ref{cor:uppersubseq} can be regarded as an evidence for Conjecture \ref{conjecture} as it has a similar form to the desired asymptote. On the other hand, due to the inexplicit choice made in the proof of Theorem \ref{thm:partialupper}, it is hard to explicitly understand from which diverging sequence $(k_j, g_j)$ we can deduce the desired asymptote. Hence it may require different approaches to make a concrete progress towards Conjecture \ref{conjecture}.

However, pseudo-Anosov mapping classes we construct in the later section (Section \ref{sec:example}) satisfy the asymptotes in Theorem \ref{thm:partialupper} and Corollary \ref{cor:uppersubseq}.

%
%

\section{Pseudo-Anosovs with specified invariant homology dimension} \label{sec:example}

To the best of the authors' knowledge, asymptotics for $L_{\C}(k, g)$ are known only when $k = 0$ (whole mapping class groups) and $k = 2g$ (Torelli groups). In this section, we construct pseudo-Anosovs $f_g \in \Mod(S_g)$ with $m(f_g) = 2g - 1$ and realizing the asymptotic for $L_{\C}(2g - 1, g)$.

From the definition of $L_{\C}(k, g)$, $L_{\C}(k, g) \le L_{\C}(k', g)$ if $k \le k'$. Since $L_{\C}(2g, g) \asymp 1 / g$ from \cite{baik2020minimal}, the lower bound in Theorem \ref{thm:main}, implies that $L_{\C}(k, g) \asymp 1 / g$ if $k$ behaves like $2g$; for instance, $k \ge 2g - C$ for some constant $C > 0$. However, $L_{\C}(0, g) \asymp 1/g^2$ by \cite{GadreTsai11}. In this regard, we ask whether there is a sort of threshold for $k$ that $L_{\C}(k, g)$ becomes strictly smaller than $1 / g$, such as $1/g^2$.

As a potential approach for this question, we think of constructing pseudo-Anosovs of specified maximal invariant homology dimensions on surfaces of large genera with small asymptotic translation lengths. In order to get pseudo-Anosov maps on surfaces of large genera, some previous results (\cite{agol2016pseudo}, \cite{baik2018upper}, \cite{baik2019asymptotic}) employ a fixed hyperbolic mapping torus and consider its monodromy obtained from a fibered cone. Since the first betti number of a mapping torus of $f$ is same as $m(f) + 1$, such monodromies in a fibered cone share the same $m$-value. 

In contrast, we come up with a finite cyclic cover of a genus 2 surface to get the desired pseudo-Anosov maps on large genera surfaces as lifts of a fixed map. From the concrete estimation on how covering maps distort the distances on curve complexes \cite{aougab2018covers}, asymptotic translation length of the lift via degree $g$ covering would be at least $1/g$, up to a constant multiple. We believe that constructing such lifts with specified maximal invariant homology dimensions would help to figure out the minimal $k = k(g)$ with $L_{\C}(k, g) \asymp 1/g$.

\begin{figure}[h]
	\centering
	\begin{tikzpicture}[scale=0.7, every node/.style={scale=0.7}]
	\draw (-1.6, 0) .. controls (-1.6, 1) and (-0.6, 0.5) .. (0, 0.5) .. controls (0.6, 0.5) and (1.6, 1) .. (1.6, 0);
	\begin{scope}[rotate=180]
		\draw (-1.6, 0) .. controls (-1.6, 1) and (-0.6, 0.5) .. (0, 0.5) .. controls (0.6, 0.5) and (1.6, 1) .. (1.6, 0);
	\end{scope}
	
	\draw (-1.1, 0.1) .. controls (-1, -0.1) and (-0.6, -0.1) .. (-0.5, 0.1);
	\draw (-1, 0) .. controls (-0.9, 0.1) and (-0.7, 0.1) .. (-0.6, 0);
	
	\draw (1.1, 0.1) .. controls (1, -0.1) and (0.6, -0.1) .. (0.5, 0.1);
	\draw (1, 0) .. controls (0.9, 0.1) and (0.7, 0.1) .. (0.6, 0);
	
	\draw[red] (1, 0) .. controls (1.1, 0.25) and (1.5, 0.25) .. (1.6, 0);
	\draw[red, dashed] (1, 0) .. controls (1.1, -0.25) and (1.5, -0.25) .. (1.6, 0);
	
	\draw[red] (1.3, 0.3) node {$\alpha$};
	
	\draw[<-] (1.8, 0) -- (2.8, 0);
	\draw (2.3, 0) node[above] {$p_{g+1}$};
	
	\begin{scope}[shift={(4.7, 0)}, rotate=45]
		
		\draw (-2, 0) .. controls (-2, 1) and (-1, 0.5) .. (-0.75, 0.5) .. controls (-0.5, 0.5) .. (-0.5, 0.75) .. controls (-0.5, 1) and (-1, 2) .. (0, 2);
		\draw (-1.5, 0.1) .. controls (-1.4, -0.1) and (-1, -0.1) .. (-0.9, 0.1);
		\draw (-1.4, 0) .. controls (-1.3, 0.1) and (-1.1, 0.1) .. (-1, 0);
		
		\begin{scope}[shift={(1.2, 0)}]
			\draw (-1.5, 0.1) .. controls (-1.4, -0.1) and (-1, -0.1) .. (-0.9, 0.1);
			\draw (-1.4, 0) .. controls (-1.3, 0.1) and (-1.1, 0.1) .. (-1, 0);
		\end{scope}
		
		\begin{scope}[rotate=90]
			\draw (-2, 0) .. controls (-2, 1) and (-1, 0.5) .. (-0.75, 0.5) .. controls (-0.5, 0.5) .. (-0.5, 0.75) .. controls (-0.5, 1) and (-1, 2) .. (0, 2);
			\draw (-1.5, -0.1) .. controls (-1.4, 0.1) and (-1, 0.1) .. (-0.9, -0.1);
			\draw (-1.4, 0) .. controls (-1.3, -0.1) and (-1.1, -0.1) .. (-1, 0);
		\end{scope}
		
		\begin{scope}[rotate=-90]
			\draw (-2, 0) .. controls (-2, 1) and (-1, 0.5) .. (-0.75, 0.5) .. controls (-0.5, 0.5) .. (-0.5, 0.75) .. controls (-0.5, 1) and (-1.5, 1.5) .. (-1, 2);
			\draw (-1.5, 0.1) .. controls (-1.4, -0.1) and (-1, -0.1) .. (-0.9, 0.1);
			\draw (-1.4, 0) .. controls (-1.3, 0.1) and (-1.1, 0.1) .. (-1, 0);
		\end{scope}
		\begin{scope}[rotate=180]
		\draw (-2, 1) .. controls (-1.5, 1.5) and (-1, 0.5) .. (-0.75, 0.5) .. controls (-0.5, 0.5) .. (-0.5, 0.75) .. controls (-0.5, 1) and (-1, 2) .. (0, 2);
		\end{scope}
		
		\draw[red, dashed] (-0.53, 0.53) .. controls (-0.53, 0.33) and (-0.27, 0.07) .. (-0.07, 0.07);
		\draw[red] (-0.53, 0.53) .. controls (-0.33, 0.53) and (-0.07, 0.27) .. (-0.07, 0.07);
		
		\begin{scope}[rotate=-90]
			\draw[red] (-0.53, 0.53) .. controls (-0.53, 0.33) and (-0.27, 0.07) .. (-0.07, 0.07);
			\draw[red, dashed] (-0.53, 0.53) .. controls (-0.33, 0.53) and (-0.07, 0.27) .. (-0.07, 0.07);
		\end{scope}
		
		\begin{scope}[rotate=90]
		\draw[red, dashed] (-0.53, 0.53) .. controls (-0.53, 0.33) and (-0.27, 0.07) .. (-0.05, 0.05);
		\draw[red] (-0.53, 0.53) .. controls (-0.33, 0.53) and (-0.07, 0.27) .. (-0.05, 0.05);
		\end{scope}
		
		\begin{scope}[rotate=180]
		\draw[red] (-0.53, 0.53) .. controls (-0.53, 0.33) and (-0.27, 0.07) .. (-0.05, 0.05);
		\draw[red, dashed] (-0.53, 0.53) .. controls (-0.33, 0.53) and (-0.07, 0.27) .. (-0.05, 0.05);
		\end{scope}
		
		\filldraw (1.75, 0.25) circle(1pt);
		\filldraw (1.8, 0) circle(1pt);
		\filldraw (1.75, -0.25) circle(1pt);
	\end{scope}

	\end{tikzpicture}
	
	\caption{$g$-fold finite cyclic covering} \label{fig:covering}
\end{figure}
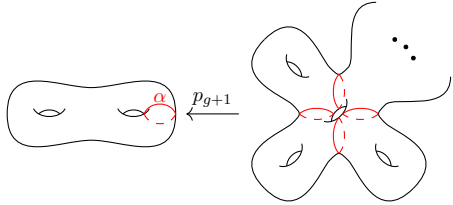

We start with a non-separating curve $\alpha$ on the genus 2 surface $S_2$, and take $g$ copies of $S_2 \setminus \alpha$ for $g > 1$. Gluing two different copies of $S_2 \setminus \alpha$ along one boundary component in a cyclic way, we obtain the finite cyclic cover $p_{g + 1}$ of degree $g$ as in Figure \ref{fig:covering}. Let us denote the resulting cover by $S_{g+1}$ since it is of genus $g+1$.

This cover $p_{g+1}$ corresponds to the kernel of the composed map  $$\pi_{1}(S_2) \xrightarrow{\hat{i}(\cdot, \alpha)} \Z \xrightarrow{\mod g} \Z / g\Z$$ where $\hat{i}(\cdot, \cdot)$ stands for the algebraic intersection number. To see this, one can observe that an element of $\pi_1(S_2)$ can be lifted to $\pi_1(S_{g+1})$ via $p_{g+1}$ if and only if its lift departs one copy of $S_2 \setminus \alpha$ and then returns to the same copy. If the lift departs and returns through the same boundary component of $S_2 \setminus \alpha$, then the element of $\pi_1(S_2)$ has the algebraic intersection number 0 with $\alpha$. Otherwise, if the lift departs and returns through different boundary components, then the algebraic intersection number is an integer multiple of $g$.

In \cite{baik2020minimal}, the first author and Shin directly constructed pseudo-Anosovs on $S_g$ that are Torelli and of small asymptotic translation lengths on curve complexes. In the following, we construct pseudo-Anosovs with specific maximal invariant homology dimensions and satisfying the upper bound provided in Theorem \ref{thm:partialupper} and Corollary \ref{cor:uppersubseq}. As a result, we obtain Theorem \ref{thm:2g+1}. Our strategy is to fix a suitable pseudo-Anosov on $S_2$ and then to lift it via $p_{g+1}$. Due to the symmetry of the covering, we can find a number of invariant homology classes proportional to the degree of the cover.

\subsection*{Proof of Theorem \ref{thm:2g+1}}

	The last assertion is a direct consequence of the first assertion and Theorem \ref{thm:main}. By Theorem \ref{thm:main}, there exists $C' > 0$ so that $L_{\C}(2g - 1, g) \ge {C' \over g}$ for all $g > 1$. Hence, it remains to show the existence of the desired pseudo-Anosovs.
	
\begin{figure}[h]
		\centering
		\begin{tikzpicture}[scale=0.7, every node/.style={scale=0.7}]
		\draw (-1.6, 0) .. controls (-1.6, 1) and (-0.6, 0.5) .. (0, 0.5) .. controls (0.6, 0.5) and (1.6, 1) .. (1.6, 0);
		\begin{scope}[rotate=180]
			\draw (-1.6, 0) .. controls (-1.6, 1) and (-0.6, 0.5) .. (0, 0.5) .. controls (0.6, 0.5) and (1.6, 1) .. (1.6, 0);
		\end{scope}
		
		\draw (-1.1, 0.1) .. controls (-1, -0.1) and (-0.6, -0.1) .. (-0.5, 0.1);
		\draw (-1, 0) .. controls (-0.9, 0.1) and (-0.7, 0.1) .. (-0.6, 0);
		
		\draw (1.1, 0.1) .. controls (1, -0.1) and (0.6, -0.1) .. (0.5, 0.1);
		\draw (1, 0) .. controls (0.9, 0.1) and (0.7, 0.1) .. (0.6, 0);
		
		\draw[red] (1, 0) .. controls (1.1, 0.25) and (1.5, 0.25) .. (1.6, 0);
		\draw[red, dashed] (1, 0) .. controls (1.1, -0.25) and (1.5, -0.25) .. (1.6, 0);
		
		\draw[red] (1.3, 0.3) node {$\alpha$};
		
		\draw[blue] (0, 0.5) .. controls (0.2, 0.3) and (0.2, -0.3) .. (0, -0.5);
		\draw[blue, dashed] (0, 0.5) .. controls (-0.2, 0.3) and (-0.2, -0.3) .. (0, -0.5);
		
		\draw[blue] (0, -0.5) node[below] {$\beta$};

		\draw[<-] (1.8, 0) -- (2.8, 0);
		\draw (2.3, 0) node[above] {$p_{g+1}$};
		
		\begin{scope}[shift={(4.7, 0)}, rotate=45]
		
			\draw (-2, 0) .. controls (-2, 1) and (-1, 0.5) .. (-0.75, 0.5) .. controls (-0.5, 0.5) .. (-0.5, 0.75) .. controls (-0.5, 1) and (-1, 2) .. (0, 2);
			\draw (-1.5, 0.1) .. controls (-1.4, -0.1) and (-1, -0.1) .. (-0.9, 0.1);
			\draw (-1.4, 0) .. controls (-1.3, 0.1) and (-1.1, 0.1) .. (-1, 0);
		
			\begin{scope}[shift={(1.2, 0)}]
				\draw (-1.5, 0.1) .. controls (-1.4, -0.1) and (-1, -0.1) .. (-0.9, 0.1);
				\draw (-1.4, 0) .. controls (-1.3, 0.1) and (-1.1, 0.1) .. (-1, 0);
			\end{scope}
		
			\begin{scope}[rotate=90]
				\draw (-2, 0) .. controls (-2, 1) and (-1, 0.5) .. (-0.75, 0.5) .. controls (-0.5, 0.5) .. (-0.5, 0.75) .. controls (-0.5, 1) and (-1, 2) .. (0, 2);
				\draw (-1.5, -0.1) .. controls (-1.4, 0.1) and (-1, 0.1) .. (-0.9, -0.1);
				\draw (-1.4, 0) .. controls (-1.3, -0.1) and (-1.1, -0.1) .. (-1, 0);
			\end{scope}
		
			\begin{scope}[rotate=-90]
				\draw (-2, 0) .. controls (-2, 1) and (-1, 0.5) .. (-0.75, 0.5) .. controls (-0.5, 0.5) .. (-0.5, 0.75) .. controls (-0.5, 1) and (-1.5, 1.5) .. (-1, 2);
				\draw (-1.5, 0.1) .. controls (-1.4, -0.1) and (-1, -0.1) .. (-0.9, 0.1);
				\draw (-1.4, 0) .. controls (-1.3, 0.1) and (-1.1, 0.1) .. (-1, 0);
			\end{scope}
			\begin{scope}[rotate=180]
				\draw (-2, 1) .. controls (-1.5, 1.5) and (-1, 0.5) .. (-0.75, 0.5) .. controls (-0.5, 0.5) .. (-0.5, 0.75) .. controls (-0.5, 1) and (-1, 2) .. (0, 2);
			\end{scope}
			
			\draw[red, dashed] (-0.53, 0.53) .. controls (-0.53, 0.33) and (-0.27, 0.07) .. (-0.07, 0.07);
			\draw[red] (-0.53, 0.53) .. controls (-0.33, 0.53) and (-0.07, 0.27) .. (-0.07, 0.07);
		
			\begin{scope}[rotate=-90]
				\draw[red] (-0.53, 0.53) .. controls (-0.53, 0.33) and (-0.27, 0.07) .. (-0.07, 0.07);
				\draw[red, dashed] (-0.53, 0.53) .. controls (-0.33, 0.53) and (-0.07, 0.27) .. (-0.07, 0.07);
			\end{scope}
		
			\begin{scope}[rotate=90]
				\draw[red, dashed] (-0.53, 0.53) .. controls (-0.53, 0.33) and (-0.27, 0.07) .. (-0.05, 0.05);
				\draw[red] (-0.53, 0.53) .. controls (-0.33, 0.53) and (-0.07, 0.27) .. (-0.05, 0.05);
			\end{scope}
		
			\begin{scope}[rotate=180]
				\draw[red] (-0.53, 0.53) .. controls (-0.53, 0.33) and (-0.27, 0.07) .. (-0.05, 0.05);
				\draw[red, dashed] (-0.53, 0.53) .. controls (-0.33, 0.53) and (-0.07, 0.27) .. (-0.05, 0.05);
			\end{scope}
			
			\draw[blue] (-0.68, 0.5) .. controls (-0.48, 0.3) and (-0.48, -0.3) .. (-0.68, -0.5);
			\draw[blue, dashed] (-0.68, 0.5) .. controls (-0.88, 0.3) and (-0.88, -0.3) .. (-0.68, -0.5);
			
			\begin{scope}[rotate=90]
				\draw[blue] (-0.68, 0.5) .. controls (-0.48, 0.3) and (-0.48, -0.3) .. (-0.68, -0.5);
				\draw[blue, dashed] (-0.68, 0.5) .. controls (-0.88, 0.3) and (-0.88, -0.3) .. (-0.68, -0.5);
			\end{scope}
			
			\begin{scope}[rotate=-90]
			\draw[blue, dashed] (-0.68, 0.5) .. controls (-0.48, 0.3) and (-0.48, -0.3) .. (-0.68, -0.5);
			\draw[blue] (-0.68, 0.5) .. controls (-0.88, 0.3) and (-0.88, -0.3) .. (-0.68, -0.5);
			\end{scope}
		
			\filldraw (1.75, 0.25) circle(1pt);
			\filldraw (1.8, 0) circle(1pt);
			\filldraw (1.75, -0.25) circle(1pt);
		\end{scope}
		
		\end{tikzpicture}
		
		\caption{A non-separating curve $\alpha$ and a separating curve $\beta$ on $S_2$ with $\alpha \cap \beta = \emptyset$} \label{fig:alphabeta}
	\end{figure}
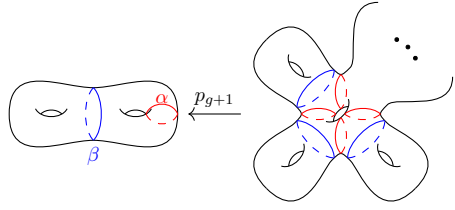

	Fix $g > 1$ and in the rest of the proof, we simply denote by $p := p_{g+1}$. Let $\alpha$ be a non-separating curve on $S_2$ and let $\beta$ be a separating curve on $S_2$ disjoint from $\alpha$. Then each lift of $\beta$ through $p$ is also separating. For instance, see Figure \ref{fig:alphabeta}.
	
	Now let $\varphi$ be a Torelli pseudo-Anosov mapping class in $\Mod(S_2)$ such that $d_{\C}(\beta, \varphi \beta) \ge 3$. Note that this means that $\beta$ and $\varphi \beta$ are separating simple closed curves on $S_2$ that fill the surface. Such a $\varphi$ may be obtained as follows: let $\psi \in \Mod(S_2)$ be a pseudo-Anosov mappingn class. Since a pseudo-Anosov mapping class acts on $\C(S_2)$ as a loxodromic isometry by \cite{MasurMinsky99}, $\lim_{n \to \infty} d_{\C}(\beta, \psi \beta)/n = \ell_{\C}(\psi) > 0$. Hence replacing $\psi$ with $\psi^n$ for some large $n$ if necessary, we may assume that $d_{\C}(\beta, \psi \beta) \ge 3$. We then consider a mapping class $T_{\beta}T_{\psi \beta}^{-1}$ which is pseudo-Anosov by Thurston \cite{thurston1988geometry} or Penner \cite{penner1988construction}. Again, since $T_{\beta}T_{\psi \beta}^{-1}$ acts on $\C(S_2)$ as a loxodromic isometry, we can set $\varphi$ to be some power of $T_{\beta}T_{\psi \beta}^{-1}$ so that $\varphi$ is a pseudo-Anosov mapping class with the property that $d_{\C}(\beta, \varphi \beta) \ge 3$.
	
	Now let $f = T_{\beta}T_{\varphi \beta}^{-1}T_{\varphi \alpha}^{-1}$. Since $\beta$, $\varphi \beta$, and $\varphi \alpha$ fill the surface and $\varphi \beta \cap~\varphi \alpha = \emptyset$, $f$ is pseudo-Anosov again by Thurston \cite{thurston1988geometry} and Penner \cite{penner1988construction}. Since $\beta$ and $\varphi \beta$ are separating, $T_{\beta}$ and $T_{\varphi \beta}^{-1}$ are Torelli, in particular, they preserve the homology class $[\alpha]$ of $\alpha$. Furthermore, since $\varphi$ is Torelli, $[\varphi \alpha] = [\alpha]$, which implies that $T_{\varphi \alpha}^{-1}$ also preserves $[\alpha]$. Hence, $f$ preserves $[\alpha]$, and thus $\hat{i}(f(\cdot), \alpha) = \hat{i}(\cdot, \alpha)$. In particular, $f$ preserves the kernel of $\pi_{1}(S_2) \xrightarrow{\hat{i}(\cdot, \alpha)} \Z \xrightarrow{\mod g} \Z / g\Z$. Consequently, $f$ can be lifted through $p$.
	
	Let $\tilde{f} = T_{p^{-1}(\beta)}T_{p^{-1}(\varphi\beta)}^{-1}T_{p^{-1}(\varphi\alpha)}^{-1}$ be a lift of $f$ via $p$. We now estimate $\ell_{\C}(\tilde{f})$. Our strategy to obtain the desired upper bound for $\ell_{\C}{\tilde{F}}$ is to find a simple closed curve $\tilde{\alpha}$ such that $\tilde{\alpha}$ and its image under a sufficiently high power of $\tilde{f}$ do not fill the surface, which means that they are within distance $2$ in the curve complex. We do this by counting the number of intersections of images of $\tilde{\alpha}$ and lifts of subsurfaces.

	Recall the construction of $p$: Take $g$ copies $X_{1}, \ldots, X_{g}$ of $S_2 \setminus \alpha$ and glue $X_{i}$ and $X_{i+1}$ along one of their boundary components. Throughout, we write each index $i$ modulo $g$. Let $\tilde{\alpha} = \partial X_0 \cap \partial X_1$. That is, $\tilde{\alpha}$ be a boundary component of $X_0$ and $X_1$ where they are glued. Due to the construction, $\tilde{\alpha}$ is a lift of $\alpha$.

\begin{figure}[h]
		\centering
		\begin{tikzpicture}[scale=0.7, every node/.style={scale=0.7}]
		\draw (-1.6, 0) .. controls (-1.6, 1) and (-0.6, 0.5) .. (0, 0.5) .. controls (0.6, 0.5) and (1.6, 1) .. (1.6, 0);
		\begin{scope}[rotate=180]
		\draw (-1.6, 0) .. controls (-1.6, 1) and (-0.6, 0.5) .. (0, 0.5) .. controls (0.6, 0.5) and (1.6, 1) .. (1.6, 0);
		\end{scope}
		
		\draw (-1.1, 0.1) .. controls (-1, -0.1) and (-0.6, -0.1) .. (-0.5, 0.1);
		\draw (-1, 0) .. controls (-0.9, 0.1) and (-0.7, 0.1) .. (-0.6, 0);
		
		\draw (1.1, 0.1) .. controls (1, -0.1) and (0.6, -0.1) .. (0.5, 0.1);
		\draw (1, 0) .. controls (0.9, 0.1) and (0.7, 0.1) .. (0.6, 0);
		
		\draw[red] (1, 0) .. controls (1.1, 0.25) and (1.5, 0.25) .. (1.6, 0);
		\draw[red, dashed] (1, 0) .. controls (1.1, -0.25) and (1.5, -0.25) .. (1.6, 0);
		
		\draw[red] (1.3, 0.3) node {$\alpha$};
		
		\draw[olive, thick] (-0.7, -0.05) .. controls (-0.5, -0.4) and (1.2, -0.25).. (1.2, 0.15) .. controls (1.2, 0.3) and (0.6, 0.3) .. (0.6, 0);
		\draw[olive, thick] (-0.6, 0) .. controls (-0.6, 0.4) and (1.5, 0.8) .. (1.5, 0.15) .. controls (1.5, -0.6) and (-0.5, -0.2).. (-0.7, -0.6);
		
		\draw[olive, dashed, thick] (-0.7, -0.6) .. controls (-0.85, -0.4) and (-0.85, -0.25) .. (-0.7, -0.05);
		\draw[olive, dashed, thick] (-0.6, 0) .. controls (-0.4, -0.2) and (0.4, -0.2) .. (0.6, 0);
		
		\draw[<-] (1.8, 0) -- (2.8, 0);
		\draw (2.3, 0) node[above] {$p_{g+1}$};
		
		\begin{scope}[shift={(4.7, 0)}, rotate=45]
		
		\draw (-2, 0) .. controls (-2, 1) and (-1, 0.5) .. (-0.75, 0.5) .. controls (-0.5, 0.5) .. (-0.5, 0.75) .. controls (-0.5, 1) and (-1, 2) .. (0, 2);
		\draw (-1.5, 0.1) .. controls (-1.4, -0.1) and (-1, -0.1) .. (-0.9, 0.1);
		\draw (-1.4, 0) .. controls (-1.3, 0.1) and (-1.1, 0.1) .. (-1, 0);
		
		\begin{scope}[shift={(1.2, 0)}]
		\draw (-1.5, 0.1) .. controls (-1.4, -0.1) and (-1, -0.1) .. (-0.9, 0.1);
		\draw (-1.4, 0) .. controls (-1.3, 0.1) and (-1.1, 0.1) .. (-1, 0);
		\end{scope}
		
		\begin{scope}[rotate=90]
		\draw (-2, 0) .. controls (-2, 1) and (-1, 0.5) .. (-0.75, 0.5) .. controls (-0.5, 0.5) .. (-0.5, 0.75) .. controls (-0.5, 1) and (-1, 2) .. (0, 2);
		\draw (-1.5, -0.1) .. controls (-1.4, 0.1) and (-1, 0.1) .. (-0.9, -0.1);
		\draw (-1.4, 0) .. controls (-1.3, -0.1) and (-1.1, -0.1) .. (-1, 0);
		\end{scope}
		
		\begin{scope}[rotate=-90]
		\draw (-2, 0) .. controls (-2, 1) and (-1, 0.5) .. (-0.75, 0.5) .. controls (-0.5, 0.5) .. (-0.5, 0.75) .. controls (-0.5, 1) and (-1.5, 1.5) .. (-1, 2);
		\draw (-1.5, 0.1) .. controls (-1.4, -0.1) and (-1, -0.1) .. (-0.9, 0.1);
		\draw (-1.4, 0) .. controls (-1.3, 0.1) and (-1.1, 0.1) .. (-1, 0);
		\end{scope}
		\begin{scope}[rotate=180]
		\draw (-2, 1) .. controls (-1.5, 1.5) and (-1, 0.5) .. (-0.75, 0.5) .. controls (-0.5, 0.5) .. (-0.5, 0.75) .. controls (-0.5, 1) and (-1, 2) .. (0, 2);
		\end{scope}
		
		\draw[red, dashed] (-0.53, 0.53) .. controls (-0.53, 0.33) and (-0.27, 0.07) .. (-0.07, 0.07);
		\draw[red] (-0.53, 0.53) .. controls (-0.33, 0.53) and (-0.07, 0.27) .. (-0.07, 0.07);
		
		\begin{scope}[rotate=-90]
		\draw[red] (-0.53, 0.53) .. controls (-0.53, 0.33) and (-0.27, 0.07) .. (-0.07, 0.07);
		\draw[red, dashed] (-0.53, 0.53) .. controls (-0.33, 0.53) and (-0.07, 0.27) .. (-0.07, 0.07);
		\end{scope}
		
		\begin{scope}[rotate=90]
		\draw[red, dashed] (-0.53, 0.53) .. controls (-0.53, 0.33) and (-0.27, 0.07) .. (-0.05, 0.05);
		\draw[red] (-0.53, 0.53) .. controls (-0.33, 0.53) and (-0.07, 0.27) .. (-0.05, 0.05);
		\end{scope}
		
		\begin{scope}[rotate=180]
		\draw[red] (-0.53, 0.53) .. controls (-0.53, 0.33) and (-0.27, 0.07) .. (-0.05, 0.05);
		\draw[red, dashed] (-0.53, 0.53) .. controls (-0.33, 0.53) and (-0.07, 0.27) .. (-0.05, 0.05);
		\end{scope}
		
		\filldraw (1.75, 0.25) circle(1pt);
		\filldraw (1.8, 0) circle(1pt);
		\filldraw (1.75, -0.25) circle(1pt);
		
	
		\draw[olive, dashed, thick] (-0.7-0.5, -0.6) .. controls (-0.85-0.5, -0.4) and (-0.85-0.5, -0.25) .. (-0.7-0.5, -0.05);
		\draw[olive, thick] (-1.2, -0.05) .. controls (-1, -0.3) and (-0.2, -0.3) .. (0, -0.05);
		\draw[olive, thick] (-1.2, -0.6) .. controls (-1, -0.3) and (0, -0.3) .. (0, -1);
		\draw[olive, dashed, thick] (0, -0.05) .. controls (0.1, -0.25) and (0.1, -0.8) .. (0, -1);
		
		\begin{scope}[rotate=180]
		\draw[olive] (-1.2, -0.05) .. controls (-1, -0.3) and (-0.2, -0.3) .. (0, -0.07);
		\draw[olive] (-1.2, -0.6) .. controls (-1, -0.3) and (0, -0.3) .. (0, -1);
		\draw[olive, dashed] (0, -0.07) .. controls (0.1, -0.25) and (0.1, -0.8) .. (0, -1);
		\end{scope}
		
		\begin{scope}[rotate=-90]
		\draw[olive, dashed] (-0.7-0.5, -0.6) .. controls (-0.85-0.5, -0.4) and (-0.85-0.5, -0.25) .. (-0.7-0.5, -0.05);
		\draw[olive] (-1.2, -0.05) .. controls (-1, -0.3) and (-0.2, -0.3) .. (0, -0.2);
		\draw[olive] (-1.2, -0.6) .. controls (-1, -0.3) and (0, -0.3) .. (0, -1);
		\draw[olive, dashed] (0, -0.2) .. controls (0.1, -0.25) and (0.1, -0.8) .. (0, -1);
		\end{scope}
		
		\begin{scope}[rotate=90]
		\draw[olive, dashed] (-0.7-0.5, -0.6) .. controls (-0.85-0.5, -0.4) and (-0.85-0.5, -0.25) .. (-0.7-0.5, -0.05);
		\draw[olive] (-1.2, -0.05) .. controls (-1, -0.3) and (-0.2, -0.3) .. (0, -0.2);
		\draw[olive] (-1.2, -0.6) .. controls (-1, -0.3) and (0, -0.3) .. (-0.5, -1);
		\draw[olive, dashed] (0, -0.2) .. controls (-0.1, -0.25) and (-0.1, -0.8) .. (-0.3, -1);
		\end{scope}
		
		\end{scope}
		
		\end{tikzpicture}
		
		\caption{A curve intersecting a lift of the new curve (non-indexed one) in some $X_j$ spreads into $X_{j-1} \cup X_{j+1}$ by twisting along the lift. It describes how the image of $\tilde{\alpha}$ under multitwists is trapped in the certain number of lifts of a subsurface, as in \eqref{eqn:trapped}.} \label{fig:spreading}
	\end{figure}
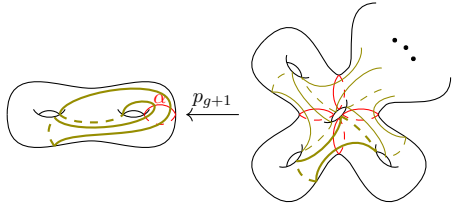

	Noting that $\hat{i}(\varphi \alpha, \alpha) = 0$ since $\varphi$ is Torelli, we get 
	\begin{equation}  \label{eqn:trapped}
		T_{p^{-1}(\varphi \alpha)}^{-1}\tilde{\alpha} \subseteq \bigcup_{j = - i(\varphi \alpha, \alpha)/2}^{i(\varphi \alpha, \alpha)/2} X_j
	\end{equation}
	 where $i(\cdot, \cdot)$ is the geometric intersection number. (cf. Figure \ref{fig:spreading}) Similarly, $\hat{i}(\varphi \beta, \alpha) = 0$ and $$T_{p^{-1}(\varphi \beta)}^{-1}T_{p^{-1}(\varphi \alpha)}^{-1}\tilde{\alpha} \subseteq \bigcup_{j = - {i(\varphi \beta, \alpha) + i(\varphi \alpha, \alpha) \over 2}}^{i(\varphi \beta, \alpha) + i(\varphi \alpha, \alpha) \over 2} X_j.$$ Since $T_{p^{-1}(\beta)}$ fixes each $X_j$, we have 
	$$
		\tilde{f}\tilde{\alpha} \subseteq \bigcup_{j = - {i(\varphi \beta, \alpha) + i(\varphi \alpha, \alpha) \over 2}}^{i(\varphi \beta, \alpha) + i(\varphi \alpha, \alpha) \over 2} X_j.
	$$
	
	Conducting this procedure inductively, we finally get $$\tilde{f}^n\tilde{\alpha} \subseteq \bigcup_{j = - n \cdot {i(\varphi \beta, \alpha) + i(\varphi \alpha, \alpha) \over 2}}^{n \cdot {i(\varphi \beta, \alpha) + i(\varphi \alpha, \alpha) \over 2}} X_j.$$ Hence, there exists $\tilde{j}$ such that $\tilde{f}^{\left\lfloor {g - 2 \over i(\varphi \beta, \alpha) + i(\varphi \alpha, \alpha)}\right\rfloor} \tilde{\alpha} \cap X_{\tilde{j}} = \emptyset$. Since there exists an essential simple closed curve in $X_{\tilde{j}}$ which is a 2-holed torus, we have $d_{\C}\left(\tilde{\alpha}, \tilde{f}^{\left\lfloor {g - 2 \over i(\varphi \beta, \alpha) + i(\varphi \alpha, \alpha)}\right\rfloor} \tilde{\alpha}\right) \le 2$ so $\ell_{\C}\left(\tilde{f}^{\left\lfloor {g - 2 \over i(\varphi \beta, \alpha) + i(\varphi \alpha, \alpha)}\right\rfloor} \right) \le 2$. This implies the following estimate. Note that $\varphi$, $i(\varphi\beta, \alpha)$, and $i(\varphi\alpha, \alpha)$ are universal quantities independent on $p$ and $g$. \begin{equation} \label{eqn:concretebound}
	\ell_{\C}(\tilde{f}) \le {2 \over \left\lfloor {g - 2 \over i(\varphi \beta, \alpha) + i(\varphi \alpha, \alpha)}\right\rfloor}\end{equation}

\begin{figure}[h]
		\centering
		\begin{tikzpicture}[scale=1, every node/.style={scale=0.7}]
			\draw (-1.6, 0) .. controls (-1.6, 1) and (-0.6, 0.5) .. (0, 0.5) .. controls (0.6, 0.5) and (1.6, 1) .. (1.6, 0);
			\begin{scope}[rotate=180]
				\draw (-1.6, 0) .. controls (-1.6, 1) and (-0.6, 0.5) .. (0, 0.5) .. controls (0.6, 0.5) and (1.6, 1) .. (1.6, 0);
			\end{scope}
		
			\draw (-1.1, 0.1) .. controls (-1, -0.1) and (-0.6, -0.1) .. (-0.5, 0.1);
			\draw (-1, 0) .. controls (-0.9, 0.1) and (-0.7, 0.1) .. (-0.6, 0);
		
			\draw (1.1, 0.1) .. controls (1, -0.1) and (0.6, -0.1) .. (0.5, 0.1);
			\draw (1, 0) .. controls (0.9, 0.1) and (0.7, 0.1) .. (0.6, 0);
		
			\draw[red, thick] (1, 0) .. controls (1.1, 0.25) and (1.5, 0.25) .. (1.6, 0);
			\draw[red, dashed, thick] (1, 0) .. controls (1.1, -0.25) and (1.5, -0.25) .. (1.6, 0);

			\draw[red] (1.3, 0.3) node {$\alpha$};
			
			\draw[teal, thick] (-1, 0) .. controls (-1.1, 0.25) and (-1.5, 0.25) .. (-1.6, 0);
			\draw[teal, dashed, thick] (-1, 0) .. controls (-1.1, -0.25) and (-1.5, -0.25) .. (-1.6, 0);
			\draw[teal] (-1.6, 0) node[left] {$\gamma$};
			
			\draw[teal, thick] (-1.3, 0) .. controls (-1.3, 0.5) and (-0.25, 0.5) .. (-0.25, 0) .. controls (-0.25, -0.4) and (-1.3, -0.4) .. (-1.3, 0);
			\draw[teal] (-0.8, -0.45) node {$\delta$};
			
			\draw[teal, thick] (1.3, 0) .. controls (1.3, 0.5) and (0.25, 0.5) .. (0.25, 0) .. controls (0.25, -0.4) and (1.3, -0.4) .. (1.3, 0);
			\draw[teal] (0.8, -0.45) node {$\eta$};
		
			\draw[blue] (0, 0.5) .. controls (0.2, 0.3) and (0.2, -0.3) .. (0, -0.5);
			\draw[blue, dashed] (0, 0.5) .. controls (-0.2, 0.3) and (-0.2, -0.3) .. (0, -0.5);
		
			\draw[blue] (0, -0.5) node[below] {$\beta$};

			\draw[<-] (1.8, 0) -- (2.8, 0);
			\draw (2.3, 0) node[above] {$p_{g+1}$};
		
			\begin{scope}[shift={(4.7, 0)}, rotate=45]
		
				\draw (-2, 0) .. controls (-2, 1) and (-1, 0.5) .. (-0.75, 0.5) .. controls (-0.5, 0.5) .. (-0.5, 0.75) .. controls (-0.5, 1) and (-1, 2) .. (0, 2);
				\draw (-1.5, 0.1) .. controls (-1.4, -0.1) and (-1, -0.1) .. (-0.9, 0.1);
				\draw (-1.4, 0) .. controls (-1.3, 0.1) and (-1.1, 0.1) .. (-1, 0);
		
				\begin{scope}[shift={(1.2, 0)}]
					\draw (-1.5, 0.1) .. controls (-1.4, -0.1) and (-1, -0.1) .. (-0.9, 0.1);
					\draw (-1.4, 0) .. controls (-1.3, 0.1) and (-1.1, 0.1) .. (-1, 0);
				\end{scope}
		
				\begin{scope}[rotate=90]
					\draw (-2, 0) .. controls (-2, 1) and (-1, 0.5) .. (-0.75, 0.5) .. controls (-0.5, 0.5) .. (-0.5, 0.75) .. controls (-0.5, 1) and (-1, 2) .. (0, 2);
					\draw (-1.5, -0.1) .. controls (-1.4, 0.1) and (-1, 0.1) .. (-0.9, -0.1);
					\draw (-1.4, 0) .. controls (-1.3, -0.1) and (-1.1, -0.1) .. (-1, 0);
				\end{scope}
		
				\begin{scope}[rotate=-90]
					\draw (-2, 0) .. controls (-2, 1) and (-1, 0.5) .. (-0.75, 0.5) .. controls (-0.5, 0.5) .. (-0.5, 0.75) .. controls (-0.5, 1) and (-1.5, 1.5) .. (-1, 2);
					\draw (-1.5, 0.1) .. controls (-1.4, -0.1) and (-1, -0.1) .. (-0.9, 0.1);
					\draw (-1.4, 0) .. controls (-1.3, 0.1) and (-1.1, 0.1) .. (-1, 0);
				\end{scope}
				\begin{scope}[rotate=180]
					\draw (-2, 1) .. controls (-1.5, 1.5) and (-1, 0.5) .. (-0.75, 0.5) .. controls (-0.5, 0.5) .. (-0.5, 0.75) .. controls (-0.5, 1) and (-1, 2) .. (0, 2);
				\end{scope}
		
				\draw[red, dashed, thick] (-0.53, 0.53) .. controls (-0.53, 0.33) and (-0.27, 0.07) .. (-0.07, 0.07);
				\draw[red, thick] (-0.53, 0.53) .. controls (-0.33, 0.53) and (-0.07, 0.27) .. (-0.07, 0.07);
				\draw[red] (-0.55, 0.55) node[left] {$\tilde{\alpha}$};
		
				\draw[teal, thick] (-1.4, 0) .. controls (-1.5, 0.25) and (-1.9, 0.25) .. (-2, 0);
				\draw[teal, dashed, thick] (-1.4, 0) .. controls (-1.5, -0.25) and (-1.9, -0.25) .. (-2, 0);
				\draw[teal, thick] (-1.7, 0) .. controls (-1.7, 0.5) and (-0.65, 0.5) .. (-0.65, 0) .. controls (-0.65, -0.4) and (-1.7, -0.4) .. (-1.7, 0);
				\draw[teal] (-2, 0) node[left] {$\tilde{\gamma}_1$};
				\draw[teal] (-1.5, -0.45) node {$\tilde{\delta}_1$};
				
				\begin{scope}[rotate=90]
					\draw[teal, dashed, thick] (-1.4, 0) .. controls (-1.5, 0.25) and (-1.9, 0.25) .. (-2, 0);
					\draw[teal, thick] (-1.4, 0) .. controls (-1.5, -0.25) and (-1.9, -0.25) .. (-2, 0);
					\draw[teal, thick] (-1.7, 0) .. controls (-1.7, 0.4) and (-0.65, 0.4) .. (-0.65, 0) .. controls (-0.65, -0.5) and (-1.7, -0.5) .. (-1.7, 0);
					\draw[teal] (-2, 0) node[right] {$\tilde{\gamma}_2$};
					\draw[teal] (-1.6, 0.45) node {$\tilde{\delta}_2$};
				\end{scope}
				
				\begin{scope}[rotate=-90]
					\draw[teal, thick] (-1.4, 0) .. controls (-1.5, 0.25) and (-1.9, 0.25) .. (-2, 0);
					\draw[teal, dashed, thick] (-1.4, 0) .. controls (-1.5, -0.25) and (-1.9, -0.25) .. (-2, 0);
					\draw[teal, thick] (-1.7, 0) .. controls (-1.7, 0.5) and (-0.65, 0.5) .. (-0.65, 0) .. controls (-0.65, -0.4) and (-1.7, -0.4) .. (-1.7, 0);
					\draw[teal] (-2, 0) node[left] {$\tilde{\gamma}_0$};
					\draw[teal] (-1.5, -0.45) node {$\tilde{\delta}_0$};
				\end{scope}
				
				\draw[teal, thick] (-1.7+1.175, 0) .. controls (-1.7+1.175, 0.5) and (-0.65+1.175, 0.5) .. (-0.65+1.175, 0) .. controls (-0.65+1.175, -0.4) and (-1.7+1.175, -0.4) .. (-1.7+1.175, 0);
				\draw[teal] (0.565, 0) node[right] {$\tilde{\eta}$};
		
				\filldraw (1.75, 0.25) circle(1pt);
				\filldraw (1.8, 0) circle(1pt);
				\filldraw (1.75, -0.25) circle(1pt);
			\end{scope}
		
		\end{tikzpicture}
		
		\caption{Basis for $H_1(S_{g+1})$} \label{fig:basis}
	\end{figure}
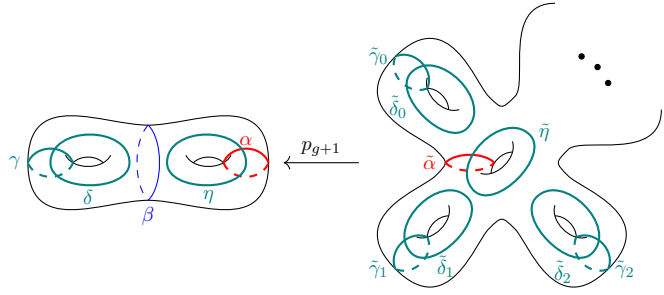

	We now show that $m(\tilde{f}) = 2g + 1$. Recall that $\beta \subseteq S_2$ is a separating curve and $\alpha \subseteq S_2$ is a non-separating curve disjoint from $\beta$. Temporarily let us denote $Y$ be a component of $S_2 \setminus \beta$ that does not contain $\alpha$. Let $\gamma$ and $\delta$ be non-separating simple closed curves whose homology classes form a basis for $H_1(Y) \cong \R^2$. Let us also denote $\eta$ be a non-separating curve on $S_2 \setminus \beta$ with $i(\eta, \alpha)=1$. Then $p^{-1}(\gamma)$, $p^{-1}(\delta)$, $p^{-1}(\eta)$, and one component of $p^{-1}(\alpha)$ form a basis for $H_1(S_{g+1}) \cong \R^{2g+2}$. See Figure \ref{fig:basis} for instance.
	
	Let $\tilde{\gamma}_j = p^{-1}(\gamma) \cap X_j$, $\tilde{\delta}_j = p^{-1}(\delta) \cap X_j$, and $\tilde{\eta} = p^{-1}(\eta)$. Further, let $\tilde{\alpha} = \partial X_0 \cap \partial X_1$ which is a component of $p^{-1}(\alpha)$. Since $\varphi$ is Torelli, it has a lift $\tilde{\varphi}$ through $p$. Hence, homology classes $\left\lbrace \left[\tilde{\varphi}\tilde{\gamma}_j\right], \left[\tilde{\varphi}\tilde{\delta}_j\right], \left[\tilde{\varphi}\tilde{\eta}\right], \left[\tilde{\varphi}\tilde{\alpha}\right] \right\rbrace$ also form a basis for $H_1(S_{g+1})$.
	
	Recall that $\tilde{f} = T_{p^{-1}(\beta)}T_{p^{-1}(\varphi\beta)}^{-1}T_{p^{-1}(\varphi\alpha)}^{-1}$. Since $\gamma \cap (\alpha \cup \beta) = \emptyset$, we have $\tilde{\varphi}\tilde{\gamma}_j \cap (p^{-1}(\varphi\alpha) \cap p^{-1}(\varphi\beta)) = \emptyset$. Here, note that $\tilde{\varphi}\tilde{\gamma}_j$ is a lift of $\varphi\gamma$ which is a component of $p^{-1}(\varphi \gamma)$. Hence it follows that $\left[\tilde{f}\tilde{\varphi}\tilde{\gamma}_j\right] = \left[T_{p^{-1}(\beta)} \tilde{\varphi}\tilde{\gamma}_j\right]$. Since each component of $p^{-1}(\beta)$, which is a lift of $\beta$, is separating, $T_{p^{-1}(\beta)}$ is Torelli. As a result, $\left[\tilde{f}\tilde{\varphi}\tilde{\gamma}_j\right] = \left[\tilde{\varphi}\tilde{\gamma}_j\right]$. Similarly, we have $\left[\tilde{f}\tilde{\varphi}\tilde{\delta}_j\right] = \left[\tilde{\varphi}\tilde{\delta}_j\right]$.
	
	Now we consider $\left[ \tilde{f}\tilde{\varphi}\tilde{\alpha}\right]$. As $\tilde{\varphi}\tilde{\alpha}$ is a lift of $\varphi\alpha$, $T_{p^{-1}(\varphi\alpha)}^{-1}\tilde{\varphi}\tilde{\alpha} =\tilde{\varphi}\tilde{\alpha}$. Furthermore, since $\alpha \cap \beta = \emptyset$, $\tilde{\varphi}\tilde{\alpha}$, a lift of $\varphi \alpha$, does not intersect $p^{-1}(\varphi \beta)$. It implies that $T_{p^{-1}(\varphi\beta)}^{-1}\tilde{\varphi}\tilde{\alpha} = \tilde{\varphi}\tilde{\alpha}$. Finally, since $T_{p^{-1}(\beta)}$ is Torelli again, we conclude $\left[\tilde{f}\tilde{\varphi}\tilde{\alpha}\right] = \left[\tilde{\varphi}\tilde{\alpha}\right]$.
	
	So far, we have proved $m(\tilde{f}) \ge 2g + 1$. Suppose to the contrary that $m(\tilde{f}) = 2g + 2$, which means that $\tilde{f}$ is Torelli. Then $\tilde{f}\tilde{\varphi}\tilde{\eta}$ should be homologous to $\tilde{\varphi}\tilde{\eta}$. It implies $\left[T_{p^{-1}(\varphi\beta)}^{-1}T_{p^{-1}(\varphi\alpha)}^{-1}\tilde{\varphi}\tilde{\eta}\right] = \left[\tilde{\varphi}\tilde{\eta}\right]$ since $T_{p^{-1}(\beta)}$ is Torelli. As any two components of $p^{-1}(\varphi \alpha)$ bound a subsurface, they are homologous. In particular, as $\tilde{\varphi}\tilde{\alpha}$ is a component of $p^{-1}(\varphi \alpha)$, each of its components is homologous to $\tilde{\varphi}\tilde{\alpha}$. Hence, $\left[T_{p^{-1}(\varphi\beta)}^{-1}T_{p^{-1}(\varphi\alpha)}^{-1}\tilde{\varphi}\tilde{\eta}\right] = \left[T_{p^{-1}(\varphi\beta)}^{-1}T_{\tilde{\varphi}\tilde{\alpha}}^{-g}\tilde{\varphi}\tilde{\eta}\right]$. Noting that $T_{\tilde{\alpha}}^{-g}\tilde{\eta}$ can be isotoped into arbitrary neighborhood of $\tilde{\alpha} \cup \tilde{\eta}$, $T_{\tilde{\varphi}\tilde{\alpha}}^{-g}\tilde{\varphi}\tilde{\eta}$ can also be isotoped into arbitrary neighborhood of $\tilde{\varphi}\tilde{\alpha} \cup \tilde{\varphi}\tilde{\eta}$. Since $\tilde{\varphi}\tilde{\alpha} \cup \tilde{\varphi}\tilde{\eta}$ and $p^{-1}(\varphi\beta)$ are disjoint compact sets, we have $T_{p^{-1}(\varphi\beta)}^{-1}T_{\tilde{\varphi}\tilde{\alpha}}^{-g}\tilde{\varphi}\tilde{\eta} = T_{\tilde{\varphi}\tilde{\alpha}}^{-g}\tilde{\varphi}\tilde{\eta}$. Summing up the above argument, we obtain $$
	\left[\tilde{\varphi}\tilde{\eta}\right]  = \left[\tilde{f}\tilde{\varphi}\tilde{\eta}\right]
	= \left[T_{p^{-1}(\varphi\beta)}^{-1}T_{p^{-1}(\varphi\alpha)}^{-1}\tilde{\varphi}\tilde{\eta}\right]
	= \left[T_{p^{-1}(\varphi\beta)}^{-1}T_{\tilde{\varphi}\tilde{\alpha}}^{-g}\tilde{\varphi}\tilde{\eta}\right]
	= \left[T_{\tilde{\varphi}\tilde{\alpha}}^{-g}\tilde{\varphi}\tilde{\eta}\right]$$ where the first equality is the assumption. However, $$\left[T_{\tilde{\varphi}\tilde{\alpha}}^{-g}\tilde{\varphi}\tilde{\eta}\right] = \left[ \tilde{\varphi}\tilde{\eta}\right] - g \cdot \hat{i}(\tilde{\varphi}\tilde{\eta}, \tilde{\varphi} \tilde{\alpha}) \left[\tilde{\varphi}\tilde{\alpha}\right]$$ which implies that $\hat{i}(\tilde{\varphi}\tilde{\eta}, \tilde{\varphi}\tilde{\alpha}) = 0$. It contradicts our choice of $\eta$ that $i(\tilde{\eta}, \tilde{\alpha}) =~1$. Therefore, $m(\tilde{f}) = 2g + 1$. Setting $f_{g+1} = \tilde{f}$ completes the proof of Theorem \ref{thm:2g+1}.

	\medskip
	
The lower bound on $\ell_{\C}(f_g)$ for $f_g$ constructed in the proof can also be calculated in a concrete way by Aougab-Patel-Taylor \cite{aougab2018covers} as follows. $${\ell_{\C}(f) \over (g-1) \cdot 80\cdot 2^{13}e^{54}\pi} \le \ell_{\C}(f_{g})$$

\begin{remark} \label{rmk:torelli}

In the proof, all figures describe one specific example. Any choice of $\alpha$, $\beta$, $\gamma$, $\delta$, and $\eta$ works if it satisfies the condition we provide. That is, 
\begin{itemize}
	\item $\alpha$ and $\beta$ are non-separating and separating curves on $S_2$ respectively, and are disjoint;
	\item $\gamma$ and $\delta$ are non-separating simple closed curves that form a basis for the first homology group of  the component of $S_2 \setminus \beta$ disjoint from $\alpha$;
	\item $\eta$ is a non-separating curve on $S_2 \setminus \beta$ with $i(\eta, \alpha) = 1$.
\end{itemize}

 Furthermore, if we modify the map on $S_2$ to be $f = T_{\beta}T_{\varphi\beta}^{-1}$, then its lift via $p_{g+1}$ is Torelli, which gives another proof of $L_{\C}(2g, g) \asymp {1 \over g}$.

\end{remark}

\section{Small translation length and Normal generation} \label{sec:normal}
	
In this section, we discuss pseudo-Anosov mapping classes with small asymptotic translation lengths and normal generation of mapping class groups. For a general group $G$ and $g \in G$, the \emph{normal closure} $\llangle g \rrangle$ of $g$ is the smallest normal subgroup of $G$ containing $g$. The normal closure can be described in various way: $$
\llangle g \rrangle = \bigcap_{g \in N \trianglelefteq G} N = \langle hgh^{-1} : h \in G \rangle
$$ In a particular case that $\llangle g \rrangle = G$, we say $g$ \emph{normally generates} $G$, and $g$ is said to be a \emph{normal generator} of $G$.

Normal generators of mapping class groups of surfaces have been studied by various authors. In \cite{long1986note}, Long asked whether there is a pseudo-Anosov that normally generates a mapping class group. This question recently answered affirmatively by Lanier and Margalit in \cite{lanier2018normal}. Indeed, they showed that there is a universal constant so that pseudo-Anosovs with stretch factors less than the constant should be normal generators. Then the asymptotic $L_{\T}(0, g) \asymp 1/g$ by Penner \cite{penner1991bounds} deduces the answer. Precisely, Lanier and Margalit proved the following.

\begin{thm}[Lanier-Margalit \cite{lanier2018normal}]
	If a pseudo-Anosov $\varphi \in \Mod(S_g)$ has the stretch factor less than $\sqrt{2}$, then $\varphi$ normally generates $\Mod(S_g)$.
\end{thm}

As the logarithm of stretch factor of a pseudo-Anosov equals to the translation length of the pseudo-Anosov on the Teichm\"uller space, Lanier and Margalit's result also means that the small translation length on the Teichm\"uller space implies the normal generation of the mapping class group. One natural question in this philosophy is whether the same holds in the circumstance of curve complexes. There are several ways to formalize this question: \begin{enumerate}
	\item Is there a universal constant $C > 0$ so that if a pseudo-Anosov $\varphi \in \Mod(S_g)$ has $\ell_{\C}(\varphi) < C / g$ then $\llangle \varphi \rrangle = \Mod(S_g)$? \label{BKSW}
	
	\item Is there a universal constant $C > 0$ so that if a non-Torelli pseudo-Anosov $\varphi \in \Mod(S_g)$ has $\ell_{\C}(\varphi) < C / g$ then $\llangle \varphi \rrangle = \Mod(S_g)$? \label{nontorelli}
\end{enumerate}
Indeed, the first and the third authors of current paper, Kin and Shin stated \eqref{BKSW} in \cite[Question 1.2]{baik2019asymptotic}.

\begin{remark}
	In the above questions, the factor $1/g$ is necessary since $L_{\C}(2g, g) \asymp 1/g$ \cite{baik2020minimal} and due to Proposition \ref{prop:res}. Furthermore, we separately state above two questions in order to forbid the trivial (Torelli) case in \eqref{nontorelli} and deal with the same problem.
\end{remark}

\subsection*{Proof of Theorem \ref{prop:res}}

The family of pseudo-Anosovs constructed in Theorem \ref{thm:2g+1} actually consists of non-normal generators, that is, $\llangle f_g \rrangle \neq \Mod(S_g)$. To see this, recall that $f_g = T_{p_g^{-1}(\beta)}T_{p_g^{-1}(\varphi \beta)}^{-1}T_{p_g^{-1}(\varphi \alpha)}^{-1}$. It can be rewritten as $$f_g = T_{p_g^{-1}(\beta)} \left( \tilde{\varphi}T_{p_g^{-1}(\beta)}^{-1}\tilde{\varphi}^{-1}\right) \left(\tilde{\varphi}T_{p_g^{-1}(\alpha)}^{-1}\tilde{\varphi}^{-1}\right).$$ Hence, it follows that $\llangle f_g \rrangle \le \llangle T_{p_g^{-1}(\beta)}, T_{p_g^{-1}(\alpha)} \rrangle$ where the right-hand-side means the smallest normal subgroup containing $T_{p_g^{-1}(\beta)}$ and $T_{p_g^{-1}(\alpha)}$.

Since each component of $p_g^{-1}(\beta)$ is separating, $T_{p_g^{-1}(\beta)}$ is Torelli, namely, contained in the kernel of the symplectic representation $\Mod(S_g) \to \Sp(2g, \Z)$. Moreover, any two components of $p_g^{-1}(\alpha)$ bound an essential subsurface so they are homologous, which means that $T_{p_g^{-1}(\alpha)}$ acts samely as $T_{\tilde{\alpha}}^{g-1}$ on $H_1(S_g; \Z)$. As such, $T_{p_g^{-1}(\alpha)}$ acts trivially on the $\mod$ $(g-1)$ homology $H_1(S_g, \Z / (g-1)\Z)$. Hence, we have that the symplectic representation of $T_{p_g^{-1}(\alpha)}$ is contained in the kernel of $\Sp(2g, \Z) \to \Sp(2g, \Z / (g-1)\Z)$. Consequently, the normal closure $\llangle T_{p_g^{-1}(\beta)}, T_{p_g^{-1}(\alpha)} \rrangle$ is contained in the kernel of the composition $$\Mod(S_g) \to \Sp(2g, \Z) \to \Sp(2g, \Z / (g-1)\Z)$$ which is surjective. It follows that $$\llangle f_g \rrangle \le \llangle T_{p_g^{-1}(\beta)}, T_{p_g^{-1}(\alpha)} \rrangle \neq \Mod(S_g)$$ so $f_g$ is not a normal generator as desired.

\begin{figure}[h]
	\centering
	\begin{tikzpicture}[scale=1, every node/.style={scale=1}]
	\draw (-1.6, 0) .. controls (-1.6, 1) and (-0.6, 0.5) .. (0, 0.5) .. controls (0.6, 0.5) and (1.6, 1) .. (1.6, 0);
	\begin{scope}[rotate=180]
	\draw (-1.6, 0) .. controls (-1.6, 1) and (-0.6, 0.5) .. (0, 0.5) .. controls (0.6, 0.5) and (1.6, 1) .. (1.6, 0);
	\end{scope}
	
	\draw (-1.1, 0.1) .. controls (-1, -0.1) and (-0.6, -0.1) .. (-0.5, 0.1);
	\draw (-1, 0) .. controls (-0.9, 0.1) and (-0.7, 0.1) .. (-0.6, 0);
	
	\draw (1.1, 0.1) .. controls (1, -0.1) and (0.6, -0.1) .. (0.5, 0.1);
	\draw (1, 0) .. controls (0.9, 0.1) and (0.7, 0.1) .. (0.6, 0);
	
	\draw[red] (1, 0) .. controls (1.1, 0.25) and (1.5, 0.25) .. (1.6, 0);
	\draw[red, dashed] (1, 0) .. controls (1.1, -0.25) and (1.5, -0.25) .. (1.6, 0);
	
	\draw[red] (1.6, 0) node[right] {$\alpha$};
	
	\draw[blue] (0, 0.5) .. controls (0.2, 0.3) and (0.2, -0.3) .. (0, -0.5);
	\draw[blue, dashed] (0, 0.5) .. controls (-0.2, 0.3) and (-0.2, -0.3) .. (0, -0.5);
	
	\draw[blue] (0, -0.5) node[below] {$\beta$};
	
	\draw[olive] (-0.6, 0) .. controls (-0.6, 0.4) and (1.2, 0.4) .. (1.2, 0) .. controls (1.2, -0.5) and (-1.4, -0.5) .. (-1.4, 0) .. controls (-1.4, 0.5) and (-0.7, 0.55) .. (-0.6, 0.58);
	\draw[olive, dashed] (-0.6, 0.58) .. controls (-0.3, 0.3) and (0, 0) .. (0.6, 0);
	\draw[olive] (0.6, 0) .. controls (0.6, -0.3) and (-1.2, -0.3) .. (-1.2, 0) .. controls (-1.2, 0.6) and (1.4, 0.6) .. (1.4, 0) .. controls (1.4, -0.5) and (0.7, -0.55) .. (0.6, -0.58);
	\draw[olive, dashed] (0.6, -0.58) .. controls (0.3, -0.5) and (-0.6, -0.2) .. (-0.6, 0);
	
	\draw[olive] (-0.6, 0.58) node[above] {$\xi$};
	\end{tikzpicture}
	
	\caption{$\beta$ and $\xi$ fill the surface} \label{fig:quant}
\end{figure}
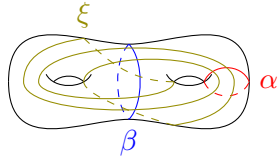

Note that we have a concrete upper bound for $\ell_{\C}(f_g)$ in (\ref{eqn:concretebound}): $$\ell_{\C}(f_g) \le {2 \over \left\lfloor {g - 3 \over i(\varphi \beta, \alpha) + i(\varphi \alpha, \alpha)}\right\rfloor} \le {2 \left( i(\varphi \beta, \alpha) + i(\varphi \alpha, \alpha)\right) \over g - 3 - \left( i(\varphi \beta, \alpha) + i(\varphi \alpha, \alpha)\right)} $$
Hence, once we fix $\alpha$, $\beta$, and $\varphi$, we get a quantitative restriction on the constant $C$ in the above questions. For instance, we can consider the configuration as in Figure \ref{fig:quant}.

Let $\lambda = T_{\xi}\beta$. As $\beta$ and $\xi$ fill the surface $S_2$, $\beta$ and $\lambda= T_{\xi}\beta$ also fill the surface. Since $\beta$ is separating, $\lambda= T_{\xi}\beta$ is also separating. Hence, due to Penner \cite{penner1988construction} or Thurston \cite{thurston1988geometry}, $\varphi=T_{\lambda}T_{\beta}^{-1}$ is a Torelli pseudo-Anosov. Furthermore, it follows that $\beta$ and $\varphi \beta$ also fill the surface. Therefore, we can construct $f_g$ as in Theorem \ref{thm:2g+1} starting with $\alpha$, $\beta$, and $\varphi$ depicted above.

Since $i(\xi, \beta) = 6$, $i(\lambda, \beta) = i(T_{\xi}\beta, \beta) = i(\xi, \beta)^2 = 36$ by \cite[Proposition 3.2]{farb2011primer}. Now from $\varphi \alpha = T_{\lambda} \alpha$ and $\varphi \beta = T_{\lambda} \beta$, we have: $$\begin{aligned}
i(\varphi \alpha, \alpha) & = i(T_{\lambda} \alpha, \alpha) = i(\lambda, \alpha)^2 & = 144 \\
i(\varphi \beta, \alpha) & = i(T_{\lambda} \beta, \alpha) = i(\lambda, \beta) i(\lambda, \alpha) & = 432
\end{aligned}$$ Hence, for the resulting $f_g$, $$\ell_{\C}(f_g) \le {1152 \over g - 579}$$ for $g > 579$. Consequently, we conclude Theorem \ref{prop:res}.

%
%

\medskip
\bibliographystyle{alpha} 

\newcommand{\etalchar}[1]{$^{#1}$}

\end{document}